\documentclass[12pt,a4paper]{amsart}
\usepackage{amssymb,amsfonts,amsthm,amsmath}
\usepackage{graphicx} 
\usepackage{color}

\theoremstyle{plain}
\newtheorem{theorem}{Theorem}

\numberwithin{equation}{section} \numberwithin{theorem}{section}
\numberwithin{lemma}{section} \numberwithin{definition}{section}
\numberwithin{corollary}{section}
\numberwithin{proposition}{section}
\begin{document}
\title{On integer solutions to $x^5 - (x+1)^5 - (x+2)^5 + (x+3)^5 = 5^m + 5^n$}
\author{Geoffrey B Campbell}
\address{Mathematical Sciences Institute \\
         The Australian National University \\
         ACT, 0200, Australia}
 \email{Geoffrey.Campbell@anu.edu.au}

\author{Aleksander Zujev}
\address{Physics Department \\
         University of California \\
         Davis, California 95616, USA}
\email{azujev@ucdavis.edu}

\begin{abstract}
We give an infinite number of integer solutions to the Diophantine equation $x^5 - (x+1)^5 - (x+2)^5 + (x+3)^5 = 5^m + 5^n$, and some solutions to some similar equations.
\end{abstract}
\maketitle

\section{A new kind of Diophantine equation with solutions}

One of us, in other investigations, found the equation
\begin{equation}
\label{eq1} 1561^5 - 1562^5 - 1563^5 + 1564^5 = 5^7 + 5^{16}.
\end{equation}
He mentioned this in a public communication on LinkedIn, which led to a discussion between the two authors of this note.

Subsequently, we found the following
\begin{theorem}For positive integers $k$,
  \begin{equation}
  \label{eq1sol}
  \left(\frac{5^{k}-3}{2}\right)^5 - \left(\frac{5^{k}-1}{2}\right)^{5} - \left(\frac{5^{k}+1}{2}\right)^5 + \left(\frac{5^{k}+3}{2}\right)^5 = 5^{k+1}+5^{3k+1}
  \end{equation}
\end{theorem}
This theorem is trivially verified in Mathematica, Maple, or an online engine such as WolframAlpha.
The equation (\ref{eq1sol}) gives us an infinite number of integer solutions to,
\begin{equation}
\label{eq2} x^5 - (x+1)^5 - (x+2)^5 + (x+3)^5 = 5^m + 5^n,
\end{equation}
which is the same as saying
\begin{equation}
\label{eq3} 10 (2 x + 3) (2x^2 + 6x + 7) = 5^m + 5^n.
\end{equation}
The set of integer cases of \ref{eq1sol}, starts with:
\begin{equation}
\label{eq4} 1^5 - 2^5 - 3^5 + 4^5 = 5^3 + 5^4,
\end{equation}
\begin{equation}
\label{eq5} 11^5 - 12^5 - 13^5 + 14^5 = 5^4 + 5^7,
\end{equation}
\begin{equation}
\label{eq6} 61^5 - 62^5 - 63^5 + 64^5 = 5^5 + 5^{10},
\end{equation}
\begin{equation}
\label{eq7} 311^5 - 312^5 - 313^5 + 314^5 = 5^6 + 5^{13},
\end{equation}
\begin{equation}
\label{eq8} 1561^5 - 1562^5 - 1563^5 + 1564^5 = 5^7 + 5^{16},
\end{equation}
\begin{equation}
\label{eq9} 7811^5 - 7812^5 - 7813^5 + 7814^5 = 5^8 + 5^{19},
\end{equation}
\begin{equation}
\label{eq10} 39061^5 - 39062^5 - 39063^5 + 39064^5 = 5^9 + 5^{22},
\end{equation}
\begin{equation}
\label{eq11} 195311^5 - 195312^5 - 195313^5 + 195314^5 = 5^{10} + 5^{25}.
\end{equation}
Evidently the above theorem is a new result. Extensive searches of the literature have not unearthed the above result, nor does it seem there are known similar types of identity where such high powers are involved in the right side. The obvious internet place to look online is the book by Tito Piezas III \cite{tP2010}. Also, Andrews \cite{gA1971} is a good modern text to browse and Dickson's very old but detailed text \cite{lD1999} is also a good resource.

\section{Further new similar results with solutions}

The authors also have found similar results to the above theorem for third and fourth power consecutive numbers. For example,
\begin{theorem}
  \begin{equation}
  \label{eq1sol2}
  \left(\frac{6^{k}-3}{2}\right)^3 - \left(\frac{6^{k}-1}{2}\right)^{3} - \left(\frac{6^{k}+1}{2}\right)^3 + \left(\frac{6^{k}+3}{2}\right)^3 = 6^{k+1}.
  \end{equation}
\end{theorem}
\begin{theorem}
  \begin{equation}
  \label{eq1sol3}
  \left(\frac{6^{k}-3}{2}\right)^4 - \left(\frac{6^{k}-1}{2}\right)^{4} - \left(\frac{6^{k}+1}{2}\right)^4 + \left(\frac{6^{k}+3}{2}\right)^4 = 6^{2k+1} + 10.
  \end{equation}
\end{theorem}
These results, as mentioned, seem to not be in the literature, and are perhaps, strangely new. Some first few cases of (\ref{eq1sol2}) and (\ref{eq1sol3}) are given by:
\begin{equation}
\label{eq12}
\left(\frac{3}{2}\right)^3 - \left(\frac{5}{2}\right)^{3} - \left(\frac{7}{2}\right)^3 + \left(\frac{9}{2}\right)^3 = 6^2,
\end{equation}
\begin{equation}
\label{eq13}
\left(\frac{33}{2}\right)^3 - \left(\frac{35}{2}\right)^{3} - \left(\frac{37}{2}\right)^3 + \left(\frac{39}{2}\right)^3 = 6^3,
\end{equation}
\begin{equation}
\label{eq14}
\left(\frac{213}{2}\right)^3 - \left(\frac{215}{2}\right)^{3} - \left(\frac{217}{2}\right)^3 + \left(\frac{219}{2}\right)^3 = 6^{4},
\end{equation}
\begin{equation}
\label{eq15}
\left(\frac{1233}{2}\right)^3 - \left(\frac{1235}{2}\right)^{3} - \left(\frac{1237}{2}\right)^3 + \left(\frac{1239}{2}\right)^3 = 6^5,
\end{equation}
\begin{equation}
\label{eq16}
\left(\frac{7773}{2}\right)^3 - \left(\frac{7775}{2}\right)^{3} - \left(\frac{7777}{2}\right)^3 + \left(\frac{7779}{2}\right)^3 = 6^6.
\end{equation}
Also,
\begin{equation}
\label{eq17}
\left(\frac{3}{2}\right)^4 - \left(\frac{5}{2}\right)^{4} - \left(\frac{7}{2}\right)^4 + \left(\frac{9}{2}\right)^4 = 6^3 + 10,
\end{equation}
\begin{equation}
\left(\frac{33}{2}\right)^4 - \left(\frac{35}{2}\right)^{4} - \left(\frac{37}{2}\right)^4 + \left(\frac{39}{2}\right)^4 = 6^5 + 10,
\end{equation}
\begin{equation}
\label{eq19}
\left(\frac{213}{2}\right)^4 - \left(\frac{215}{2}\right)^{4} - \left(\frac{217}{2}\right)^4 + \left(\frac{219}{2}\right)^4 = 6^7 + 10,
\end{equation}
\begin{equation}
\label{eq20}
\left(\frac{1233}{2}\right)^4 - \left(\frac{1235}{2}\right)^{4} - \left(\frac{1237}{2}\right)^4 + \left(\frac{1239}{2}\right)^4 = 6^9 + 10,
\end{equation}
\begin{equation}
\label{eq21}
\left(\frac{7773}{2}\right)^4 - \left(\frac{7775}{2}\right)^{4} - \left(\frac{7777}{2}\right)^4 + \left(\frac{7779}{2}\right)^4 = 6^{11} + 10.
\end{equation}
It would seem that these type of identities may be easily derived, and that the results of the present note are not the final word on this topic.


\begin{thebibliography}{99}
\bibitem{mA1972}
ABRAMOWITZ, M., and STEGUN, I. Handbook of Mathematical Functions,  Dover Publications Inc., New York, 1972.
\bibitem{gA1971}
ANDREWS, G. E.  Number Theory.  W. B. Saunders, Philadelphia,
1971.  (Reprinted: Hindustan Publishing Co., New Delhi, 1984)
\bibitem{lD1999}
DICKSON, L. E. History of the Theory of Numbers, Vol II, Ch XXII, page 644, originally published 1919 by Carnegie Inst of Washington, reprinted by The American Mathematical Society 1999.
\bibitem{tP2010}
PIEZAS III, T.  A Collection of Algebraic Identities, Email: tpiezas@gmail.com, Website: https://sites.google.com/site/tpiezas/Home, 2010.  (Cited as at: 29th February, 2016)
\end{thebibliography}
\end{document}